\documentclass[10pt]{article}
\usepackage{fullpage}
\usepackage{setspace}
\usepackage{parskip}
\usepackage{titlesec}
\usepackage[section]{placeins}
\usepackage{xcolor}
\usepackage{breakcites}
\usepackage{amsthm}
\usepackage{amsmath}
\usepackage{lineno}
\usepackage{hyphenat}
\usepackage{enumerate}

\newcommand{\sign}{  {\rm sign}  }
\PassOptionsToPackage{hyphens}{url}
\usepackage[colorlinks = true,
            linkcolor = blue,
            urlcolor  = blue,
            citecolor = blue,
            anchorcolor = blue]{hyperref}
\usepackage{etoolbox}
\makeatletter
\patchcmd\@combinedblfloats{\box\@outputbox}{\unvbox\@outputbox}{}{%
}%
\makeatother

\usepackage{natbib}

\renewenvironment{abstract}
  {{\bfseries\noindent{\abstractname}\par\nobreak}\footnotesize}
  {\bigskip}

\titlespacing{\section}{0pt}{*3}{*1}
\titlespacing{\subsection}{0pt}{*2}{*0.5}
\titlespacing{\subsubsection}{0pt}{*1.5}{0pt}

\newtheorem{thm}{Theorem}
\usepackage{authblk}

\usepackage{graphicx}
\usepackage[space]{grffile}
\usepackage{latexsym}
\usepackage{textcomp}
\usepackage{longtable}
\usepackage{tabulary}
\usepackage{booktabs,array,multirow}
\usepackage{amsfonts,amsmath,amssymb}
\providecommand\citet{\cite}
\providecommand\citep{\cite}

\newif\iflatexml\latexmlfalse

\AtBeginDocument{\DeclareGraphicsExtensions{.pdf,.PDF,.eps,.EPS,.png,.PNG,.tif,.TIF,.jpg,.JPG,.jpeg,.JPEG}}

\usepackage[utf8]{inputenc}
\usepackage[ngerman,english]{babel}

\usepackage{float}
\newcommand{\logit}{  {\rm logit}  }

\usepackage{color, colortbl}
	
\definecolor{Gray}{gray}{0.9}
\newcolumntype{g}{>{\columncolor{Gray}}c}

\begin{document}

\title{The oracle property of the generalized outcome adaptive lasso}

\author[]{Ismaila Baldé}%
\affil[]{Department of mathematics and statistics, Université de Moncton, Canada}%
\vspace{-1em}
\date{\today}
\begingroup
\let\center\flushleft
\let\endcenter\endflushleft
\maketitle
\endgroup

\selectlanguage{english}
\vspace{-0.5cm}
\begin{abstract}

The generalized outcome-adaptive lasso (GOAL) is a variable selection for high-dimensional causal inference proposed by Baldé et al. [2023, {\em Biometrics} {\bfseries 79(1)}, 514--520]. When the dimension is high, it is now well established that an ideal variable selection method should have the oracle property to ensure the optimal large sample performance. However, the oracle property of GOAL has not been proven. In this paper, we show that the GOAL estimator enjoys the oracle property. Our simulation shows that the GOAL method deals with the collinearity problem better than the oracle-like method, the outcome-adaptive lasso (OAL).

\end{abstract}%
\sloppy
\vspace{-0.5cm}
\section{Introduction}

{\label{intro}}
Let $Y$ be a continuous outcome variable, $\mathbf{X}=(X_1,\ldots, X_p)$ a potential confounders matrix and $A$ a binary treatment. Assume that all $p$ covariates $X_j$, $j=1, \ldots,p$ are measured prior to the treatment $A$ which in turn is measured prior to the outcome $Y$. We assume the propensity score (PS) model is defined as: 
$ \logit{\left\lbrace P(A=1\vert \mathbf{X})\right\rbrace}= \sum_{j=1}^p \alpha_jX_j. $
Let $\mathcal{C}$ and $\mathcal{P}$ denote indices of confounders (covariates related to both outcome and treatment) and pure predictors of outcome, respectively.
Our objective is to estimate the following PS model:
$\logit{\left\lbrace \pi ( X,\hat{\alpha})\right\rbrace}=\sum_{j\in \mathcal{C}} \hat{\alpha}_jX_j +\sum_{j\in \mathcal{P}} \hat{\alpha}_jX_j.$
The negative log-likelihood function of  $\alpha$ is given by
$\ell_n (\alpha; A,\mathbf{X})=\sum_{i=1}^n \left\lbrace-a_i(x_i^T \alpha)+\log\left( 1+e^{x_i^T \alpha} \right) \right\rbrace.$
Baldé et al. (2023) proposed the generalized outcome adaptive lasso (GOAL):
$\hat{\alpha}(GOAL) = \arg \min_{\alpha}  \left[ \ell_n(\alpha;A,\mathbf{X})  + \lambda_1 \sum_{j=1}^p \hat{w}_j|\alpha_j | + \lambda_2 \sum_{j=1}^p\alpha^2_j\right]$,
where $\hat{w}_j=\left|{\hat{\beta}_j}^{ols}\right| ^{-\gamma}$ such that $\gamma>1$ and  $(\hat{\beta}_A^{ols},\hat{\beta}^{ols})=\arg \min_{(\beta_A,\beta)}  \left\| Y-\beta_AA-\mathbf{X}\beta \right\|_2 ^2$.
GOAL is designed to improve OAL (Shortreed and Ertefaie, 2017) for high-dimensional data analysis. Baldé (2022) conjectured that GOAL must satisfy the oracle property. 
In this paper, we show that GOAL enjoys the oracle property with a proper choice of $\lambda_1$, $\lambda_2$ and $\gamma$. That is, GOAL performs as well as if the true underlying model were known in advance. The oracle property is particularly important for a variable selection method when the dimension is high to ensure optimal large sample performance (Zou, 2006; Zou and Zhang, 2009).
\vspace{-0.3cm}
\section{Statistical theory \label{stheory}}
Let $\mathcal{A}=\mathcal{C} \cup \mathcal{P}=\{ 1,2,\ldots, p_0\}$ be the indices of desirable covariates to include in the estimated PS. Let $\mathcal{A}^c=\mathcal{I} \cup \mathcal{S}=\{ p_0+1,p_0+2,\ldots, p_0+(p-p_0)\}$ be the indices of  of covariates to exclude, where $\mathcal{I}$ and $\mathcal{S}$ are the indices of pure predictors of treatment and spurious covariates (covariates that are unrelated to both outcome and treatment), respectively. We write the Fisher information (FI) matrix as 
\vspace{-0.2cm}
\small{
\begin{equation}
\mathbf{F}(\alpha^*)=E\left( \phi^{''}(x^T\alpha^*)xx^T\right)=
\begin{pmatrix}
\mathbf{F}_{11} & \mathbf{F}_{12}\\
\mathbf{F}_{21} & \mathbf{F}_{22}
\end{pmatrix},
\label{fism}
\hspace{0.2cm} \mbox{$\mathbf{F}_{11}$ is the FI matrix ($p_0\times p_0$) for the parsimonious PS.}
\end{equation}}
\vspace{-0.5cm}

\begin{thm}
Assume the following regularity conditions:
\vspace{-0.2cm}
\begin{enumerate}[\hspace{0.0cm} (C.1)]
 \item The Fisher information matrix $\mathbf{F}(\alpha^*)$ defined in Equation \ref{fism} is finite and positive definite.
\item  For each $\alpha^* \in \mathbf{\Omega}$, there existe a function $M_1(x)$ and $M_2(x)$ such that for $\alpha$ in the neighborhood of $\alpha^*$, we have: 
$\left| \phi^{''}(x,\alpha) \right|\leq M_1(x)  \quad \mbox{and} \quad \left| \phi^{'''}(x,\alpha) \right|\leq M_2(x)$
such that 

$ \int M_1(x)dx <\infty \hspace{0.2cm} \mbox{and} \hspace{0.2cm} E(M_2(x)|x_j,x_k,x_l)<\infty,  \hspace{0.1cm} \forall \hspace{0.1cm} 1\leq j,k,l\leq p_0,$
where $\mathbf{\Omega}$ is an open parameter space for $\alpha$.

\item  $\lambda_1/\sqrt{n} \rightarrow 0$ and $\lambda_1n^{\gamma/2-1} \rightarrow \infty$, for $\gamma>1$.

\item $\lambda_2/\sqrt{n} \rightarrow 0$.
\end{enumerate}
Then under conditions (C.1)-(C.4) the generalized outcome-adaptive lasso estimator $\hat{\alpha}(GOAL)$ satisfies the following:
\begin{enumerate}
    \item Consistency in variable selection: $\lim_n P\{ \hat{\alpha}_j(GOAL)=0|j\in \mathcal{I} \cup \mathcal{S}\}=1$;
    \item Asymptotic normality: $\sqrt{n} \{ \hat{\alpha}(GOAL)-\alpha^*_{\mathcal{A}} \} \rightarrow_d N(0,\mathbf{F}_{11}^{-1})$.
\end{enumerate}
\label{proofproporaclegoal}
\end{thm}
\textbf{Proof of Theorem \ref{proofproporaclegoal}.}
The ideas of the proof are taken from Zou (2006), Khalili and Chen (2007), Slawski et al. (2010), and Shortreed and Ertefaie (2017). 

First, we prove the asymptotic normality. 
Let $\alpha=\alpha^*+\frac{b}{\sqrt{n}}$. Then 
define $\mathcal{G}_n(b)$ by
\begin{flalign*}
\mathcal{G}_n(b) & =\ell_n(\alpha^*+\frac{b}{\sqrt{n}};A,\mathbf{X})  + \lambda_1 \sum_{j=1}^{p_0} \hat{w}_j \left|\alpha_j^*+\frac{b}{\sqrt{n}} \right| + \lambda_2 \sum_{j=1}^{p_0} (\alpha_j^*+\frac{b}{\sqrt{n}})^2&\\
     & =\sum_{i=1}^n -a_ix_i^T \left(\alpha^*+\frac{b}{\sqrt{n}}\right)+\phi\left( {x_i^T\left(\alpha^*+\frac{b}{\sqrt{n}}\right) } \right)  
     + \lambda_1 \sum_{j=1}^{p_0} \hat{w}_j \left|\alpha_j^*+\frac{b_j}{\sqrt{n}} \right| 
     + \lambda_2 \sum_{j=1}^{p_0} \left(\alpha_j^*+\frac{b_j}{\sqrt{n}}\right)^2.
\end{flalign*}
For $b=0$, we have: 
$\mathcal{G}_n(0)  =\sum_{i=1}^n -a_i(x_i^T \alpha^*)+\phi\left( {x_i^T\alpha^*} \right)  
     + \lambda_1 \sum_{j=1}^{p_0} \hat{w}_j \left|\alpha_j^* \right| + \lambda_2 \sum_{j=1}^{p_0} {\alpha_j^*}^2.$

Define $\widehat{\mathcal{W}}(b)=\mathcal{G}_n(b)-\mathcal{G}_n(0)$. Thus
\begin{flalign*}
\widehat{\mathcal{W}}(b) & =  \sum_{i=1}^n \left\lbrace-a_ix_i^T \left(\alpha^*+\frac{b}{\sqrt{n}}\right)+\phi\left( {x_i^T\left(\alpha^*+\frac{b}{\sqrt{n}}\right) } \right) - \left( -a_i(x_i^T \alpha^*)+\phi\left( {x_i^T\alpha^*} \right) \right) \right\rbrace&\\
    & + \lambda_1  \sum_{j=1}^{p_0} \hat{w}_j \left( \left|\alpha_j^*+\frac{b_j}{\sqrt{n}} \right|- \left|\alpha_j^* \right|  \right) 
    + \lambda_2 \sum_{j=1}^{p_0}\left( \left( \alpha_j^*+\frac{b_j}{\sqrt{n}}\right)^2 -  {\alpha_j^*}^2\right) &\\
    & =  \sum_{i=1}^n \left\lbrace-a_ix_i^T \frac{b}{\sqrt{n}}+\phi\left( {x_i^T\left(\alpha^*+\frac{b}{\sqrt{n}}\right) } \right) -  \phi\left( {x_i^T\alpha^*} \right) \right\rbrace  &\\
   & + \lambda_1  \sum_{j=1}^{p_0} \hat{w}_j \left( \left|\alpha_j^*+\frac{b_j}{\sqrt{n}} \right|- \left|\alpha_j^* \right|  \right) 
    + \lambda_2 \sum_{j=1}^{p_0}\left( \left( \alpha_j^*+\frac{b_j}{\sqrt{n}}\right)^2 -  {\alpha_j^*}^2\right) &\\
    & =\mathcal{L}_n(b)+ \lambda_1  \sum_{j=1}^{p_0} \hat{w}_j \left( \left|\alpha_j^*+\frac{b_j}{\sqrt{n}} \right|- \left|\alpha_j^* \right|  \right) 
    + \lambda_2 \sum_{j=1}^{p_0}\left( \alpha_j^*\left(\frac{2b_j}{\sqrt{n}}\right) + \frac{b_j^2}{n}\right),
\end{flalign*}
where $\mathcal{L}_n(b)=\sum_{i=1}^n \left\lbrace-a_ix_i^T \frac{b}{\sqrt{n}}+\phi\left( {x_i^T\left(\alpha^*+\frac{b}{\sqrt{n}}\right) } \right) -  \phi\left( {x_i^T\alpha^*} \right) \right\rbrace$. 

By using the second-order Taylor expansion of $\mathcal{L}_n(b)$ around $b=0$, we have:
\begin{flalign*}
  \mathcal{L}_n(b)&=\sum_{i=1}^n \left(-a_ix_i^T \frac{b}{\sqrt{n}} \right)+  \sum_{i=1}^n \left(\phi^{'}\left( {x_i^T\alpha^*} \right)   
  \left(\frac{x_i^T b}{\sqrt{n}} \right)\right) + \frac{1}{2} \sum_{i=1}^n \phi^{''}\left( {x_i^T\alpha^*}   \right) 
  b^T\left(\frac{x_i x_i^T}{n} \right)b 
  + n^{-\frac{3}{2}} \sum_{i=1}^n \frac{1}{6} \phi^{'''}\left( {x_i^T\alpha^*}\right) \left(  x_i^Tb \right)^3 &\\
  &= -\sum_{i=1}^n \left(a_i -\phi^{'}\left( {x_i^T\alpha^*} \right)  \right) \frac{x_i^T b}{\sqrt{n}}+ \frac{1}{2} \sum_{i=1}^n \phi^{''}\left( {x_i^T\alpha^*}   \right) 
  b^T\left(\frac{x_i x_i^T}{n} \right)b 
  + n^{-\frac{3}{2}} \sum_{i=1}^n \frac{1}{6} \phi^{'''}\left( {x_i^T\alpha^*}\right) \left(  x_i^Tb \right)^3. 
\end{flalign*}

Thus, we can rewrite $\widehat{\mathcal{W}}(b)$ as
$\widehat{\mathcal{W}}(b)=R_1^{(n)}+R_2^{(n)}+R_3^{(n)}+R_4^{(n)}+R_5^{(n)},$
with 
$$R_1^{(n)}=-\sum_{i=1}^n \left(a_i -\phi^{'}\left( {x_i^T\alpha^*} \right)  \right) \frac{x_i^T b}{\sqrt{n}}, \quad    
R_2^{(n)}=\frac{1}{2} \sum_{i=1}^n \phi^{''}\left( {x_i^T\alpha^*}   \right) 
  b^T\left(\frac{x_i x_i^T}{n} \right)b, \quad R_3^{(n)}=n^{-\frac{3}{2}} \sum_{i=1}^n \frac{1}{6} \phi^{'''}\left( {x_i^T\alpha^*}\right) \left(  x_i^Tb \right)^3, $$
$$R_4^{(n)}=\lambda_1  \sum_{j=1}^{p_0} \hat{w}_j \left( \left|\alpha_j^*+\frac{b_j}{\sqrt{n}} \right|- \left|\alpha_j^* \right|  \right), \quad
R_5^{(n)}= \lambda_2 \sum_{j=1}^{p_0}\left( \alpha_j^*\left(\frac{2b_j}{\sqrt{n}}\right) + \frac{b_j^2}{n}\right).$$
By applying the central limit theorem and laws of large numbers, we have:
$$R_1^{(n)}=-\sum_{i=1}^n \left(a_i -\phi^{'}\left( {x_i^T\alpha^*} \right)  \right) \frac{x_i^T b}{\sqrt{n}} \rightarrow_d b^TZ, \quad Z \sim N(0, \mathbf{F}(\alpha^*)).$$
For the term $R_2^{(n)}$, we observe that 
$$\sum_{i=1}^n \phi^{''}\left( {x_i^T\alpha^*}   \right) 
  \left(\frac{x_i x_i^T}{n} \right) \rightarrow_p \mathbf{F}(\alpha^*), \quad \mbox{thus} \quad
R_2^{(n)}=\frac{1}{2} \sum_{i=1}^n \phi^{''}\left( {x_i^T\alpha^*}   \right) 
  b^T\left(\frac{x_i x_i^T}{n} \right)b \rightarrow_p \frac{1}{2} b^T\{\mathbf{F}(\alpha^*)\} b.$$
By the condition (C.2), we observe that
\begin{equation}
 6\sqrt{n} \left\lbrace R_3^{(n)} \right\rbrace=\frac{1}{n}\sum_{i=1}^n  \phi^{'''}\left( {x_i^T\Tilde{\alpha}}_*\right) \left(  x_i^Tb \right)^3  \leq \sum_{i=1}^n \frac{1}{n}M_2\left(x_i\right)|x_i^Tb|^3 \rightarrow_p  E\left(M_2(x)\left|x^Tb\right|^3\right) <\infty,
 \label{eqr3}
  \end{equation}
  where $\Tilde{\alpha}_*$ is between $\alpha^*$ and $\alpha^*+\frac{b}{\sqrt{n}}$. Equation \ref{eqr3}  show that $R_3^{(n)}$ is bounded. 
  The behavior of $R_4^{(n)}$ and $R_5^{(n)}$ depend on the covariate type. If a covariate $X_j$ is a confounder ($j\in \mathcal{C}$) or a pure predictor of the outcome ($j\in \mathcal{P}$), this is $\alpha_j^*\neq 0$ since  $j\in \mathcal{A}=\mathcal{C} \cup \mathcal{P}$. If $\alpha_j^*\neq 0$, then we have: 
$$ \frac{\lambda_1 }{\sqrt{n}}   \hat{w}_j \sqrt{n} \left( \left|\alpha_j^*+\frac{b_j}{\sqrt{n}} \right|- \left|\alpha_j^* \right|  \right)= \left( \frac{\lambda_1 }{\sqrt{n}} \right)  \left(  \hat{w}_j \right) \left( \sqrt{n} \left( \left|\alpha_j^*+\frac{b_j}{\sqrt{n}} \right|- \left|\alpha_j^* \right|  \right) \right),$$ 
with $ \frac{\lambda_1 }{\sqrt{n}} \rightarrow_p 0$,  $\hat{w}_j \rightarrow_p \left|\beta_j^*\right| ^{-\gamma}$ and $\sqrt{n} \left( \left|\alpha_j^*+\frac{b_j}{\sqrt{n}} \right|- \left|\alpha_j^* \right|  \right) 
\rightarrow_p b_j \sign(\alpha_j^*).$

By using the Slutsky's theorem we have: 
 $ \frac{\lambda_1 }{\sqrt{n}}   \hat{w}_j \sqrt{n} \left( \left|\alpha_j^*+\frac{b_j}{\sqrt{n}} \right|- \left|\alpha_j^* \right|  \right) \rightarrow_p 0.$

For the behavior of $R_5^{(n)}$, we have
 $ \left(\frac{\lambda_2 }{\sqrt{n}}\right) \left(\sqrt{n}\left( \alpha_j^*\left(\frac{2b_j}{\sqrt{n}}\right) + \frac{b_j^2}{n}\right) \right)\rightarrow_p 0,$
since $ \frac{\lambda_2 }{\sqrt{n}} \rightarrow_p 0$ by assumption and $\sqrt{n}\left( \alpha_j^*\left(\frac{2b_j}{\sqrt{n}}\right) + \frac{b_j^2}{n} \right)  \rightarrow_p 2 \alpha_j^* b_j$ and then using the Slutsky's theorem.

Using the convexity of $\widehat{\mathcal{W}}(b)$ and following the epi-convergence results of Geyer (1994), we have:
$$ \arg\min \widehat{\mathcal{W}}(b) \rightarrow_d \arg\min \mathcal{W}(b), \quad
\mbox{where} \quad \arg\min \widehat{\mathcal{W}}(b)=\sqrt{n}(\hat{\alpha}-\alpha^*).
$$
Thus, again by applying the theorem of Slutsky, we have: $$\arg\min \mathcal{W}(b_{\mathcal{A}})=\mathbf{F}_{11}^{-1}Z_{\mathcal{A}}, \quad Z_{\mathcal{A}} \sim N(0, \mathbf{F}_{11}).$$
Finally, we prove the asymptotic normality part. 

Now we show the consistency in variable selection part, i.e. 
$\lim_n P\{ \hat{\alpha}_{\mathcal{A}^c}=0 \}=1$. Let $\alpha=(\alpha_{\mathcal{A}},\alpha_{\mathcal{A}^c})$ and define a penalized negative log-likelihood function as 
\begin{flalign*}
\tilde{\ell}_n(\alpha_{\mathcal{A}},\alpha_{\mathcal{A}^c})&=  \ell_n(\alpha_{\mathcal{A}},\alpha_{\mathcal{A}^c})  + \lambda_1 \sum_{j=1}^p \hat{w}_j|\alpha_j | + \lambda_2 \sum_{j=1}^p\alpha^2_j
=  \ell_n(\alpha_{\mathcal{A}},\alpha_{\mathcal{A}^c})  + \lambda_1 \sum_{ j\in \mathcal{A}\cup \mathcal{A}^c} \hat{w}_j|\alpha_j | + \lambda_2 \sum_{  j\in \mathcal{A}\cup \mathcal{A}^c  }\alpha^2_j.
\end{flalign*}
Thus, to prove sparsity it suffices to show that $\tilde{\ell}_n(\alpha_{\mathcal{A}},\alpha_{\mathcal{A}^c})-\tilde{\ell}_n(\alpha_{\mathcal{A}},0)>0$
with probability tending to 1 as $n \rightarrow \infty$. We observe that 
\begin{flalign*}
\tilde{\ell}_n(\alpha_{\mathcal{A}},\alpha_{\mathcal{A}^c})-\tilde{\ell}_n(\alpha_{\mathcal{A}},0)&=
 \left[ \ell_n(\alpha_{\mathcal{A}},\alpha_{\mathcal{A}^c})-\ell_n(\alpha_{\mathcal{A}},0)\right]
 +\lambda_1   \left[\sum_{ j\in \mathcal{A}\cup \mathcal{A}^c} \hat{w}_j|\alpha_j |  - \sum_{ j\in \mathcal{A}} \hat{w}_j|\alpha_j | \right] + \lambda_2  \left[ \sum_{  j\in \mathcal{A}\cup \mathcal{A}^c}  \alpha^2_j   -\sum_{ j\in  \mathcal{A}}  \alpha^2_j \right]&\\
 &=\left[ \ell_n(\alpha_{\mathcal{A}},\alpha_{\mathcal{A}^c})-\ell_n(\alpha_{\mathcal{A}},0)\right]+\lambda_1 \sum_{ j\in \mathcal{A}^c} \hat{w}_j|\alpha_j | + \lambda_2 \sum_{j\in \mathcal{A}^c  }\alpha^2_j.
 \end{flalign*}
By the mean value theorem, we have:
$\ell_n(\alpha_{\mathcal{A}},\alpha_{\mathcal{A}^c})-\ell_n(\alpha_{\mathcal{A}},0)= \left[  \frac{\partial{\ell_n(\alpha_{\mathcal{A}},\xi)}}
   {\partial {\alpha_{\mathcal{A}^c}} } \right]^T 
    \alpha_{\mathcal{A}^c}, 
$
for some $\| \xi\| \leq \|\alpha_{\mathcal{A}^c} \|$.

By the mean value theorem again, we have: 
\begin{flalign*}
\left \| \frac{\partial{\ell_n(\alpha_{\mathcal{A}},\xi)}}
   {\partial {\alpha_{\mathcal{A}^c}} } -   \frac{\partial{\ell_n(\alpha^*_{\mathcal{A}},0)}}
   {\partial {\alpha_{\mathcal{A}^c}} }  \right \|&= \left\| \left[ \frac{\partial{\ell_n(\alpha_{\mathcal{A}},\xi)}}
   {\partial {\alpha_{\mathcal{A}^c}} } -   \frac{\partial{\ell_n(\alpha_{\mathcal{A}},0)}}
   {\partial {\alpha_{\mathcal{A}^c}} }  \right]  +  \left[   \frac{\partial{\ell_n(\alpha_{\mathcal{A}},0)}}
   {\partial {\alpha_{\mathcal{A}^c}} }- \frac{\partial{\ell_n(\alpha^*_{\mathcal{A}},0)}}
   {\partial {\alpha_{\mathcal{A}^c}} }
   \right]\right\| &\\
   & \leq \left\|  \frac{\partial{\ell_n(\alpha_{\mathcal{A}},\xi)}}
   {\partial {\alpha_{\mathcal{A}^c}} } -   \frac{\partial{\ell_n(\alpha_{\mathcal{A}},0)}}
   {\partial {\alpha_{\mathcal{A}^c}} } \right\|  +  \left\|   \frac{\partial{\ell_n(\alpha_{\mathcal{A}},0)}}
   {\partial {\alpha_{\mathcal{A}^c}} }- \frac{\partial{\ell_n(\alpha^*_{\mathcal{A}},0)}}
   {\partial {\alpha_{\mathcal{A}^c}} }
   \right\|&\\
   & \leq \left( \sum_{i=1}^n M_1(x_i) \right) \left\| \xi \right\| + \left( \sum_{i=1}^n M_1(x_i) \right) \left\| \alpha_{\mathcal{A}}-\alpha^*_{\mathcal{A}} \right\|, \quad \mbox{by condition (C.2)} &\\
   & =\left( \left\| \xi \right\| +  \left\| \alpha_{\mathcal{A}}-\alpha^*_{\mathcal{A}} \right\| \right) \left( \sum_{i=1}^n M_1(x_i) \right)=\left( \left\| \xi \right\| +  \left\| \alpha_{\mathcal{A}}-\alpha^*_{\mathcal{A}} \right\| \right) O_p(n).
\end{flalign*}
Thus, the limiting behavior of $\left( \left\| \xi \right\| +  \left\| \alpha_{\mathcal{A}}-\alpha^*_{\mathcal{A}} \right\| \right) O_p(n)$ depends wether $j\in \mathcal{S}$ or $j\in \mathcal{I}$. 

If $j\in \mathcal{S}$, we have $\left\| \xi \right\| \leq  \left\| \alpha_{\mathcal{S}} \right\|= O_p(n^{-1/2})$. Thus,
$$\left \| \frac{\partial{\ell_n(\alpha_{\mathcal{A}},\xi)}}
   {\partial {\alpha_{\mathcal{A}^c}} } -   \frac{\partial{\ell_n(\alpha^*_{\mathcal{A}},0)}}
   {\partial {\alpha_{\mathcal{A}^c}} }  \right \| \leq  \left( \left\| \xi \right\| +  \left\| \alpha_{\mathcal{A}}-\alpha^*_{\mathcal{A}} \right\| \right) O_p(n)= O_p(n^{1/2}).$$

If $j\in \mathcal{I}$, we have $\left\| \xi \right\| \leq  \left\| \alpha_{\mathcal{I}} \right\|=O_p(1)$. Thus,
$$ \left \| \frac{\partial{\ell_n(\alpha_{\mathcal{A}},\xi)}}
   {\partial {\alpha_{\mathcal{A}^c}} } -   \frac{\partial{\ell_n(\alpha^*_{\mathcal{A}},0)}}
   {\partial {\alpha_{\mathcal{A}^c}} }  \right \| \leq  \left( \left\| \xi \right\| +  \left\| \alpha_{\mathcal{A}}-\alpha^*_{\mathcal{A}} \right\| \right) O_p(n)=O_p(n).$$
Then, for $j\in \mathcal{A}^c =\mathcal{I} \cup \mathcal{S}$, we have $\left \| \frac{\partial{\ell_n(\alpha_{\mathcal{A}},\xi)}}
   {\partial {\alpha_{\mathcal{A}^c}} } -   \frac{\partial{\ell_n(\alpha^*_{\mathcal{A}},0)}}
   {\partial {\alpha_{\mathcal{A}^c}} }  \right \| \leq O_p(n). $
   
Hence, we have: 
$\ell_n(\alpha_{\mathcal{A}},\alpha_{\mathcal{A}^c})-\ell_n(\alpha_{\mathcal{A}},0)=O_p(n) \sum_{ j\in \mathcal{A}^c  } -|\alpha_j |.$ Thus,  
{\footnotesize{
\begin{flalign}
\tilde{\ell}_n(\alpha_{\mathcal{A}},\alpha_{\mathcal{A}^c})-\tilde{\ell}_n(\alpha_{\mathcal{A}},0)&=\sum_{ j\in \mathcal{A}^c  }  \left[-|\alpha_j | O_p(n)+\lambda_1 \hat{w}_j|\alpha_j | + \lambda_2 \alpha^2_j\right]
=\sum_{ j\in \mathcal{A}^c  }  \left[-|\alpha_j | O_p(n)+\lambda_1 \left( O_p(1)n^{\gamma/2}  \right)|\alpha_j | + \lambda_2 \alpha^2_j\right]\label{eq21}.
\end{flalign}}}
By shrinking neighborhood of $0$, $|\alpha_j | O_p(n) \leq \lambda_1 \left( O_p(1)n^{\gamma/2}  \right)|\alpha_j | + \lambda_2 \alpha^2_j$ in probability. This show that $\tilde{\ell}_n(\alpha_{\mathcal{A}},\alpha_{\mathcal{A}^c})-\tilde{\ell}_n(\alpha_{\mathcal{A}},0) >0$ with probability tending to 1 as $n \rightarrow \infty$. 
Let $(\widehat{\alpha}_{\mathcal{A}}, 0)$ be the minimizer of the penalized negative log-likelihood function $\tilde{\ell}_n(\alpha_{\mathcal{A}}, 0)$, where $\tilde{\ell}_n(\alpha_{\mathcal{A}}, 0)$ is a function of $\alpha_{\mathcal{A}}$.  Now it suffices to prove
$\tilde{\ell}_n(\alpha_{\mathcal{A}},\alpha_{\mathcal{A}^c})-\tilde{\ell}_n(\widehat{\alpha}_{\mathcal{A}},0) >0$. 
By adding and subtracting $\tilde{\ell}_n(\alpha_{\mathcal{A}},0)$, we have:
\begin{flalign}
\tilde{\ell}_n(\alpha_{\mathcal{A}},\alpha_{\mathcal{A}^c})-\tilde{\ell}_n(\widehat{\alpha}_{\mathcal{A}},0)&=
\left[\tilde{\ell}_n(\alpha_{\mathcal{A}},\alpha_{\mathcal{A}^c})-\tilde{\ell}_n(\alpha_{\mathcal{A}},0) \right]+
\left[\tilde{\ell}_n(\alpha_{\mathcal{A}},0)-\tilde{\ell}_n(\widehat{\alpha}_{\mathcal{A}},0) \right] &\\
& \geq \tilde{\ell}_n(\alpha_{\mathcal{A}},\alpha_{\mathcal{A}^c})-\tilde{\ell}_n(\alpha_{\mathcal{A}},0)\label{eq31}.
\end{flalign}
By the results in Equation \ref{eq21}, the right-hand side of the Equation \ref{eq31} is positive with probability tending to 1 as $n \rightarrow \infty$. This completes the proof of the sparsity part. 
\vspace{-0.45cm}
\section{Numerical example \label{section numerical example}}
In this section, we present simulations done to study the finite sample performance of GOAL. Unlike the simulation studies conducted in Baldé et al. (2023), and in Shortreed and Ertefaie (2017), we follow Zou and Zhang (2009) to allow the intrinsic dimension (the size of $\mathcal{A}=\mathcal{C} \cup \mathcal{P}$) to diverge with the sample size as well. This makes our numerical study more challenging than simulation studies in Baldé et al. (2023) and in Shortreed and Ertefaie (2017), which considered a situation where the intrinsic dimension $|\mathcal{A}|$ is fixed. In this numerical study, we considered three methods: GOAL, OAL and Lasso. We use these methods because OAL is an oracle-like method (Shortreed and Ertefaie, 2017) while Lasso does not have the oracle property (Zou, 2006; Zou and Zhang, 2009).
Now, we describe the simulation set up to generate the data $(\mathbf{X},A,Y)$, denoting the covariates matrix, treatment and outcome, respectively. The vector ${X_i} = (X_{i1}, X_{i2},\ldots, X_{ip})$, for $i=1, 2, \ldots, n$ is simulated from a multivariate standard Gaussian distribution with pairwise correlation $\rho$. The binary treatment $A$ is simulated from a Bernoulli distribution with $\logit \{P(A_i=1)\}=\sum_{j=1}^p \alpha_jX_{ij}$. Given $\mathbf{X}$ and $A$, the continuous outcome is simulated as $Y_i=\beta_A A_i+\sum_{j=1}^p \beta_jX_{ij}+ \epsilon_i$ where  $\epsilon_i  \sim N(0,1)$. The true ATE was  $\beta_A=0$. We considered two different correlations: independent covariates ($\rho=0$) and strongly correlated covariates ($\rho=0.5$). Let $p=p_n=\lfloor 4n^{1/2}-5 \rfloor$ for $n = 100, 200, 400$. The true coefficients are $\alpha^*=(0.6 .1_q, 0.6 .1_q, 0.1_q, 0_{p-3q})$, $|\mathcal{A}|=3q$ with $q=\lfloor p_n/9 \rfloor$, and $1_m/0_m$ denote a $m$-vector of $1$'s/$0$'s. 
We follow Zou and Zhang (2009) to set $\gamma=3$ for computing the adaptive weights in the GOAL and OAL. For estimation accuracy, we used the bias, standard error (SE) and mean squared error (MSE) to compare methods. For variable selection, we evaluated methods based on the proportion of times each variable was selected  for inclusion in the PS model, where a variable was considered selected if its estimated coefficient was greater than the tolerance $10^{-8}$ (Shortreed and Ertefaie, 2017). 
\vspace{-0.2cm}
\begin{figure}[H]			
\begin{center}				
  \hspace{0cm}        $n=100 $ \hspace{4cm}  $n=200 $ \hspace{4cm}   $ n=400$
 		 	\includegraphics[width=5.41cm,height=4.5cm]{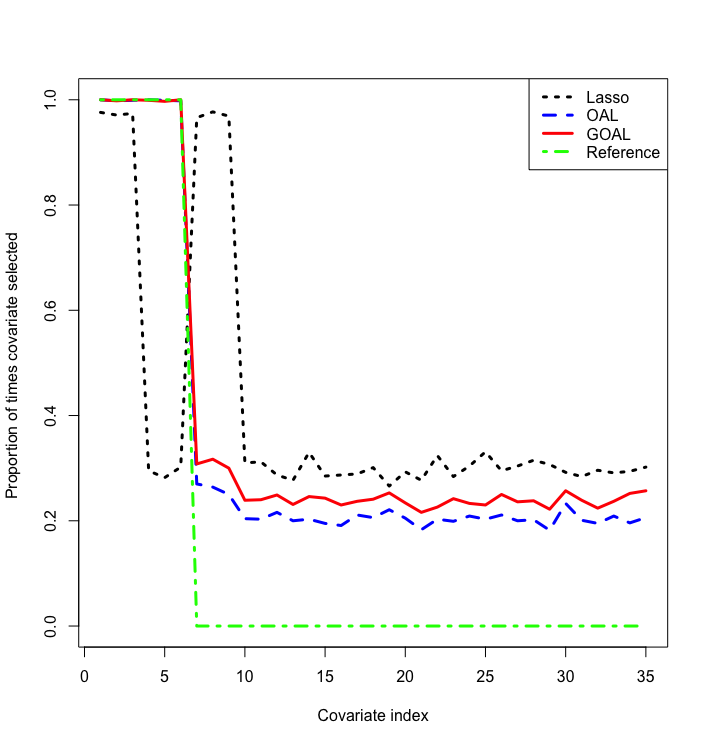}
           \includegraphics[width=5.41cm,height=4.5cm]{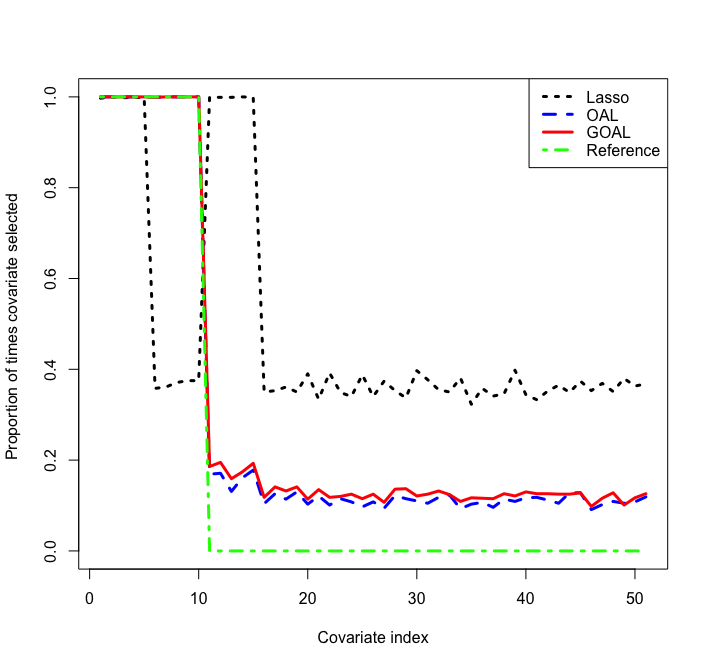}
            \includegraphics[width=5.41cm,height=4.5cm]{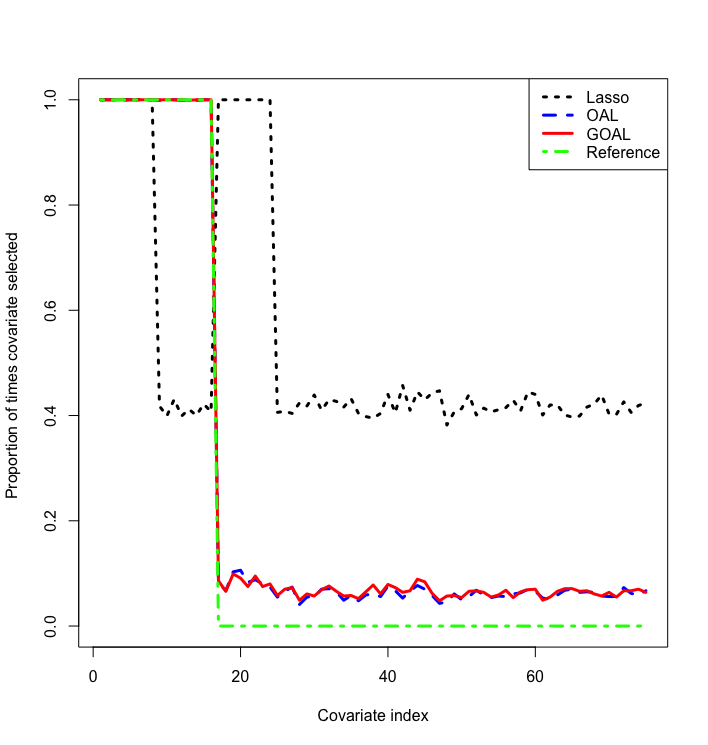}
\hspace{1cm}
\includegraphics[width=5.41cm,height=4.5cm]{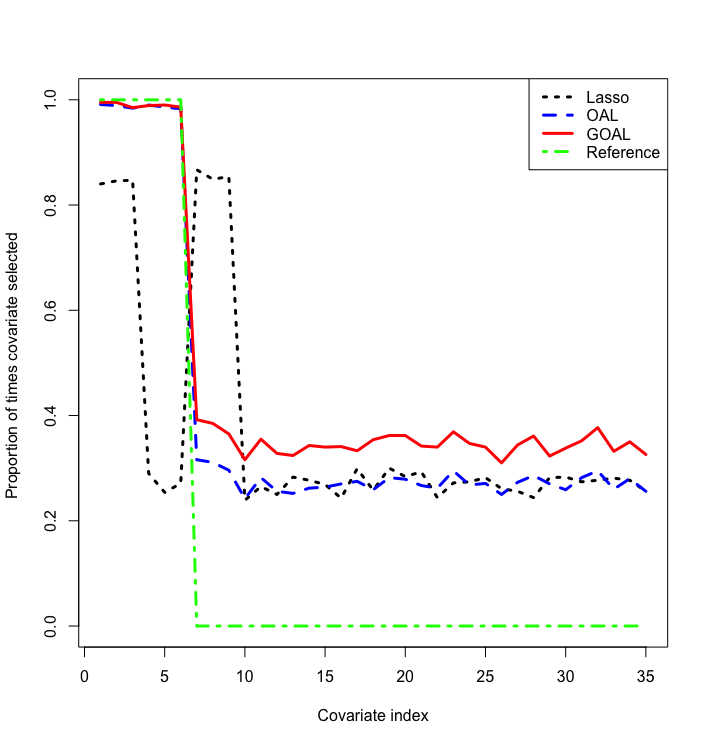}
\includegraphics[width=5.41cm,height=4.5cm]{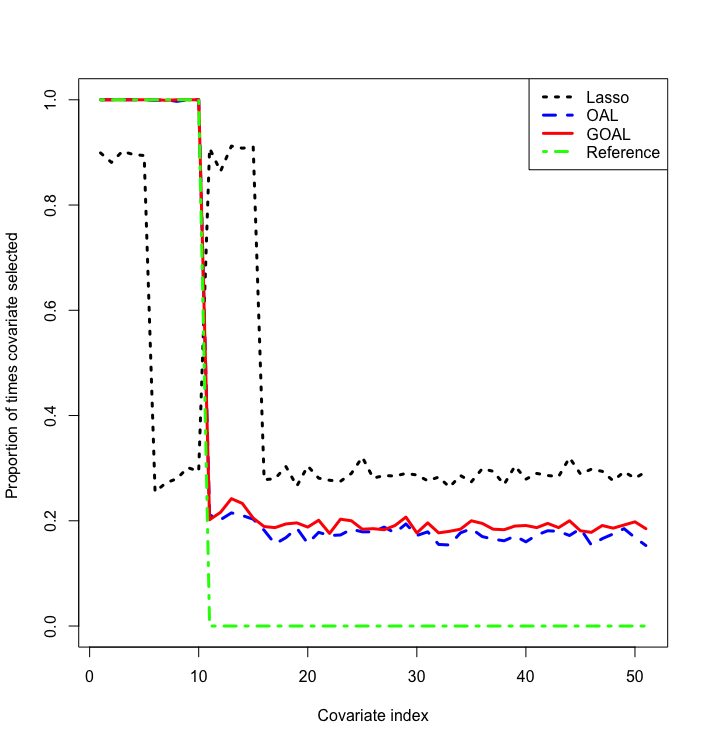}
\includegraphics[width=5.41cm,height=4.5cm]{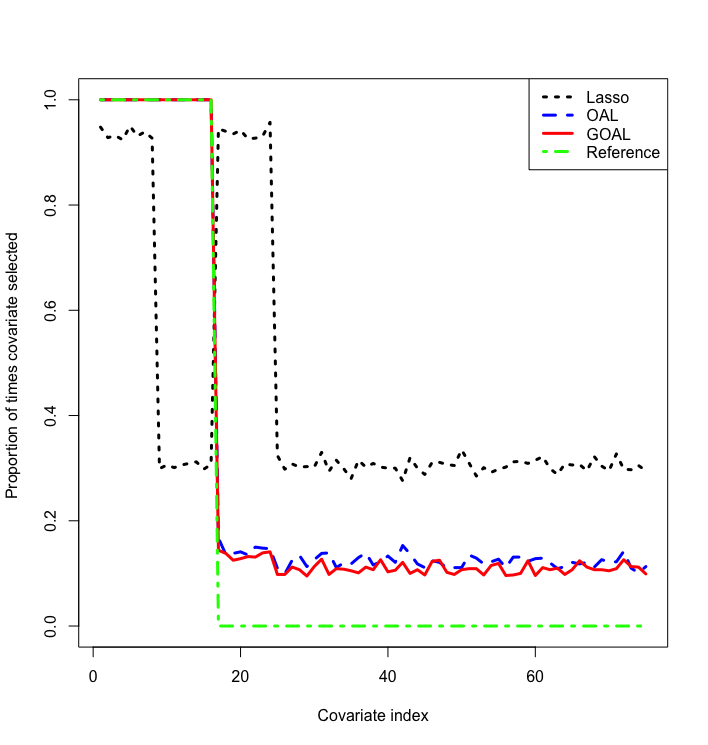}
 \vspace{-0.8cm}
\caption[]{\footnotesize{Probability of covariate being included in PS model for estimating the average treatment effect (ATE) under independent ($\rho=0$) in row 1, and strongly correlated covariates ($\rho=0.5$) in row 2}.}
\label{simres2}
\label{simres2}		
		\end{center}
		\end{figure}
\vspace{-0.9cm}
\begin{table}[H]
\footnotesize{
\begin{center}
\caption[\textbf{table}]{\footnotesize{Bias, standard error (SE) and mean squared error (MSE) of IPTW estimator for GOAL, OAL and Lasso estimator for the  average treatment effect (ATE) based on $1\, 000$ replications}.}

\begin{tabular*}{\textwidth}{@{\extracolsep\fill}|ggggggggggggggg|@{}}
\toprule
   &  &   &    & & \multicolumn{5}{@{}l}{\textbf{No correlation $\rho=0$}} & \multicolumn{5}{@{}l|}{\textbf{Strong correlation $\rho=0.5$}} \\\cmidrule{6-10}\cmidrule{11-15}
    $\mathbf{n}$    & $\mathbf{p_n}$   & $\mathbf{|\mathcal{A}|}$  & \textbf{Model} & &\textbf{Bias} && \textbf{SE} && \textbf{MSE}  &\textbf{Bias} && \textbf{SE} && \textbf{MSE} \\
  \midrule
   \rowcolor{white}
 100  & 35 & 6  & GOAL    &&  0.13  && 0.32   && 0.12   & 0.32  && 0.68   && 0.57   \\
 \rowcolor{white}
       &    &   &  OAL     &&  0.17  && 0.32   &&  0.13 &  1.40          && 0.80           && 2.60   \\
       \rowcolor{white}
      &    &   &  Lasso   &&  0.69  && 0.32   &&  0.58   &  2.93          && 0.63           && 8.97  \\
 200 & 51 & 10 & GOAL         &&  0.12   && 0.26   && 0.08  &  0.51       &&  0.84      &&   0.95  \\
     &  &  & OAL          &&  0.15   && 0.27   &&  0.10  &  2.96   && 1.14   &&  10.07   \\
     &  &  &  Lasso   &&  0.92  && 0.27   &&  0.92   &  5.51          && 0.67           && 30.79  \\
     \rowcolor{white}
   400 & 75 & 16  &GOAL         && 0.15    && 0.20   && 0.06   &  1.00    &&  1.13   && 2.27  \\
   \rowcolor{white}
     & &   & OAL          && 0.19    &&  0.22   && 0.08   &  5.56    &&  1.68   && 33.74  \\
     \rowcolor{white}
       &    &   &  Lasso   &&  1.24  && 0.24   &&  1.59   &  9.48          && 0.78           && 90.43 \\
   \bottomrule
\end{tabular*}
\label{gliores1}
\end{center}}
\end{table}
Table 1 and Figure 1 present the simulation results based on $1000$ simulations with  the $(n,p, |\mathcal{A}|)$ combinations $(100, 35, 6)$,  $(200, 51, 10)$, $(400, 75, 16)$ for both independent $(\rho=0)$ and strongly correlated $(\rho=0.5)$ covariates.

Table 1 presents the bias, standard error (SE) and mean squared error (MSE) of  GOAL, OAL and Lasso  estimators for the ATE.  In all considered scenarios, GOAL  and OAL performed much better than the Lasso method. For the independent covariates $(\rho=0)$, the MSE of GOAL and OAL decreased as the sample size increased. The method GOAL exhibited the smallest bias and MSE for all combinations $(n,p, |\mathcal{A}|, \rho)$. When covariates are strongly correlated $(\rho=0.5)$, the performance of the OAL method deteriorates with sample size. The Lasso performed the worst under the settings we considered.

Figure 1 reports the proportion of times each covariate was selected over $1000$ simulations for inclusion in the PS model when Lasso, OAL and GOAL were used to fit the PS model to estimate the ATE.  OAL and GOAL algorithms included all covariates at similar rate with high probability for confounders and pure predictors of the outcome (about 100\%), and relatively small probability for the pure predictor of treatment and spurious covariates, for all combinations $(n,p, |\mathcal{A}|, \rho)$ considered.  However, the Lasso method selects the confounders and pure predictor of treatment  with high probability and excludes  pure predictors of the outcome and spurious covariates.
\vspace{-0.2cm}
\section{Discussion}

In this paper, we studied  the statistical property of the GOAL method. We compared GOAL, OAL and Lasso using a simulation study where both the number of parameters and the intrinsic dimension $(\mathcal{A}=\mathcal{C} \cup \mathcal{P})$ diverges with the sample size. A distinctive feature of our simulation scenarios compared to many existing variable selection method for causal inference including those conducted in Baldé et al. (2023) and Shortreed and Ertefaie (2017) is that the dimension of the active set ($\mathcal{A}$) diverges with the sample size. This makes our numerical example more challenging and more appropriate for high-dimensional data analysis. GOAL and OAL outperformed the Lasso in all scenarios considered. The two oracle-like methods (GOAL and OAL) are the best when the covariates are independent ($\rho=0$) and the sample size is large (n=400). This result is expected according to the asymptotic theory for an oracle-like method (Zou and Zhang, 2009). However, OAL was less performant than GOAL when the correlation is strong ($\rho=0.5$). The GOAL method had the best performance for every combination of $(n, p, |\mathcal{A}|, \rho)$. As a result, GOAL has much better finite sample performance than the oracle-like method OAL.
\vspace{-0.4cm}
\section*{Data availability}
No data was used for the research described in the article.
\vspace{-0.4cm}
\section*{Acknowledgements}
{\label{687807}}
This work was funded by grants from  New Brunswick Innovation Foundation (NBIF).
\selectlanguage{english}
\FloatBarrier
\vspace{-0.3cm}
\section*{References}
\sloppy
\phantomsection
\label{csl:1}Baldé, I., 2022. Algorithmes de sélection de confondants en petite et grande dimensions : contextes d’application conventionnels et pour l’analyse de la médiation. Thèse. Montréal (Québec, Canada), Université du Québec à Montréal, Doctorat en mathématiques. 

\phantomsection
\label{csl:2}Baldé, I., Yang, A. Y., Lefebvre, G., 2023. {Reader Reaction to `` Outcome-adaptive lasso: Variable selection for causal inference '' by Shortreed and Ertefaie (2017)}. Biometrics 79(1), 514--520. 

\phantomsection
\label{csl:3}Khalili, A., Chen, J., 2007. {Variables selection in finite mixture of regression models}. Journal of the American Statistical Association 104, 1025--38. 

\phantomsection 
\label{csl:4}Shortreed, S. M., Ertefaie, A., 2017. {Outcome-adaptive lasso: Variable selection for causal inference}. Biometrics 73(4), 1111--1122. 

\phantomsection
\label{csl:5} Slawski, M., Castell, W. zu., Tutz, G., 2010. {Feature Selection Guided by Structural Information}. Annals of Applied Statistics,  4(2), 1056--1080. 

\phantomsection
\label{csl:6} Zou, H., 2006. {The adaptive lasso and its oracle properties}. Journal of the American Statistical Association: Series B 101, 1418--1429. 

\phantomsection
\label{csl:7} Zou, H. and Zhang, H. H. (2009).
{ On the adaptive elastic-net with a diverging number of parameters}.
The Annals of Statistics, 37, 1733--1751. 
\end{document}